\documentclass[12pt,oneside,final]{amsart}
\usepackage[cp1250]{inputenc}
\usepackage[T1]{fontenc}
\usepackage{amssymb}
\usepackage{color}
\renewcommand{\thefootnote}{}
\newtheorem{twr}{Twierdzenie}
\newtheorem{thr}[twr]{Theorem}
\newtheorem{prp}[twr]{Proposition}

\newtheorem{lm}[twr]{Lemma}

\newtheorem{crl}[twr]{Corollary}
\renewcommand{\thefootnote}{}

\newcommand{\cn}{\mathbb{C}^n}

\title[regularization of plurisubharmonic functions]{On regularization of plurisubharmonic functions near boundary points}
\author[S. Pli\'{s}]{Szymon Pli\'{s}}
\address{  Institute of Mathematics, Cracow University of Technology, Warszawska 24, 31-155
    Krak\'{o}w, Poland
}

\email{splis@pk.edu.pl}

\begin{document}\thispagestyle{empty} \footnotetext{The  author was partially supported by the NCN grant 2011/01/D/ST1/04192.
}\renewcommand{\thefootnote}{\arabic{footnote}}

\begin{abstract}
We prove in an elementary way that for  a Lipschitz domain $D\subset \cn$, all plurisubharmonic functions on $D$ can be regularized near any boundary point.
\end{abstract}

\maketitle
\section{Introduction}
 
Let $D\subset\cn$. Using a local convolution, for any plurisubharmonic function $u$ on $D$,  one can find a sequence $u_k$ of smooth plurisubharmonic functions, which decreases to $u$ on compact subsets of $D$. The purpose of this note is to show that (for regular enough $D$) it is possible to choose such a sequence near any point on the boundary.

 A domain $D\subset\cn$  will be called an $\mathcal{S}$-domain if for any plurisubharmonic function $u$ on $D$,  one can find a sequence $u_k$ of smooth plurisubharmonic functions on $D$, which decreases to $u$.\footnote{Of course this notion does not depend on coordinates so we can define it on manifolds (for example compact manifolds are $\mathcal{S}$-domains).  Theorem \ref{aproksymacjanabrzegu}  can be also formulated on manifolds.} Note that pseudoconvex domains, tube domains and Riendhard domains are $\mathcal{S}$-domains (see \cite{f-w}).

\begin{thr}\label{aproksymacjanabrzegu} Let $D\subset\cn$ be a domain with the $Lipschitz$ boundary. Then for any boundary point $P$ there is a neighbourhood $U$ of $P$ such that $D\cap U$ is  an $\mathcal{S}$-domain with $Lipschitz$ boundary.
\end{thr}

In section \ref{Example}, by a slight modification of an example from \cite{f-s}, we show the necessity of the assumption on the boundary. 

Because of the fact that not every smooth domain is an $\mathcal{S}$-domain (see \cite{f}) we have the following (surprising for the author) corollary:

\begin{crl}\label{lok}
 Being  an $\mathcal{S}$-domain is not a local property of the boundary.
\end{crl}

\section{Proof}

We need the following Lemma:
\begin{lm}\label{ciaglosc}
Let $D\subset\mathbb{R}^m$ be an open set. Let $u$ be a subharmonic function on $D$ and let $P\in D$. Let $a,b,R,C>0$, $B=\{x\in\mathbb{R}^m:|x-P|<R\}$, $K=\{(x',x_m)\in\mathbb{R}^{m-1}\times\mathbb{R}=\mathbb{R}^m:-a< x_m<-b|x'|\}$ and $B+K=\{x+y:x\in B, y\in K\}\subset D$. Assume that for any $x\in B$ and $y\in K$ $$u(x+y)\leq u(x)+\delta(|y|),$$
where $\delta:(0,+\infty)\rightarrow(0,+\infty)$ is such that $\lim_{t\rightarrow0^+}\delta(t)=0$.
Then $u$ is continuous at $P$. 
\end{lm}

\textit{Proof:} Let $(x_n)$ be any sequence in $D$ which  converges to $P$. Let $S_n=\{x\in\mathbb{R}:|x-P|=2|x_n-P|\}$, $A_n=S_n\cap(\{x_n\}+K)$ and $B_n=S_n\setminus A_n$. For $n$ large enough $x_n\in B$ and $|x_n-P|\leq \frac{a}{3} $. Hence there is a constant $\alpha>0$ (which depends only on $b$) such that  $\alpha_n=\sigma(A_n)\geq\alpha\sigma(S_n)$   where $\sigma$ is the standard measure on a sphere. Let $M_n=sup_{S_n}u$. Because $u$ is subharmonic we can estiamate:
$$u(P)\leq\sigma(S_n)^{-1}\int_{A_n}ud\sigma+\sigma(S_n)^{-1}\int_{B_n}ud\sigma$$
$$\leq \alpha_n(u(x_n)-\delta(3|x_n-P|))+(1-\alpha_n)M_n,$$
hence
$$u(x_n)\geq u(P)+\frac{1-\alpha_n}{\alpha_n}(u(P)-M_n)+\delta(3|x_n-P|).$$
Letting $n$ to $\infty$ we get
$$\varlimsup_{n\rightarrow\infty}u(x_n)\geq u(P).$$
Since $u$ is upper semicontinuous the proof is completed. $\Box$

$\newline$

The function $P_Df:=\sup\{u\in\mathcal{PSH}(D):u\leq f\}$, where $D\subset\cn$ and $f$ is a (real) function on $D$, is called a plurisubharmonic envelope of $f$. 

\textit{Proof of Theorem \ref{aproksymacjanabrzegu}:} We use the following notation $\cn\ni z=(a,x)\in(\mathbb{C}^{n-1}\times\mathbb{R})\times\mathbb{R}$. We put $B=\{a\in\mathbb{C}^{n-1}\times\mathbb{R}:|a|<1\}$. After an affine change of coordinates we can assume that there is a constant $C>0$ and a function $F:B\rightarrow[3C,4C]$ such that:\\
i) $F(a)-F(b)\leq C|a-b|$ for $a,b\in B$,\\
ii) $\partial D\cap B\times[-5C,5C]=\{(a,F(a)):a\in B\}$ and $P=(0,F(0))$,\\
iii) $0\in D$.\\
  Let $U=\{(a,x)\in\mathbb{C}^{n-1}\times\mathbb{R}:|a|^2+\left(\frac{x}{5C}\right)^2<1\}$. We will show that $\Omega=D\cap U$ is  an $\mathcal{S}$-domain. Observe that $$\partial \Omega\cap B\times[0,5C]=\{(a,\hat F(a)):a\in B\}$$ where $\hat F=\max\{F(a),5C\sqrt{1-|a|^2}\}$. By elementary calculations $\hat F(a)-\hat F(b)\leq \frac{20}{3}C|a-b|$ for $a,b\in B'=\{a\in\mathbb{C}^{n-1}\times\mathbb{R}:|a|<\frac{4}{5}\}$ and therefore the inequality $\hat F(a)-\hat F(b)\leq C'|a-b|$ holds in a neighbourhood of $\{\hat F=F\}\subset B'$ for some $C'<7C$.

Let 
$$K_\varepsilon=\{(a,x):-\varepsilon<x<-7C|a|\}.$$ Let $$\Omega_k=\Omega_k(\varepsilon)=\{z\in\Omega:z+kw\in\Omega \hbox{ for any } w\in \bar K_\varepsilon\},\hbox{ for }k=1,2.$$ Then there is a compact subset $L_\varepsilon\Subset\Omega_1$ such that  
$${\rm dist}(z,\partial\Omega)<{\rm dist}(z+w,\partial\Omega)$$ 
for $z\in \Omega_2\setminus L_\varepsilon$ and $w\in K_\varepsilon$.
 Observe that  for $\varepsilon$ small enough $\partial D\cap \Omega=\partial D\cap \Omega_2$ and the function $d=-\log({\rm dist}(\cdot,\partial\Omega))$ is plurisubharmonic in a neighbourhood of $cl_{\Omega}(\Omega\setminus\Omega_2)$. Morover, the function $d'=-\log({\rm dist}(\cdot,\partial U))$ is plurisubharmonic on the whole $\Omega$.

Let $u\in\mathcal{PSH}(\Omega)$ and let $\phi_k$ be a sequence of continuous functions on $\Omega$ which decreases to $u$. We can choose an increasing convex function $p:\mathbb{R}\rightarrow\mathbb{R}$ such that for a function $\rho=p\circ d$ we have $\lim_{z\rightarrow\partial\Omega}\rho-\phi_1=+\infty$. Put $\tilde{\phi}_k=\max\{\phi_k,\rho-k\}$.  Observe that  functions $\hat{\phi}_k=P_\Omega\tilde{\phi}_k$ are plurisubharmonic and they decrease to $u$. 

Fix $k$ and $\varepsilon>0$ such that $$cl_\Omega(\Omega\setminus\Omega_2)\subset int\{\tilde{\phi}_k=\rho'-k\},$$ where $\rho'=p\circ d'$. On $\Omega_2$ we have $$\tilde\phi_k\geq P_{\Omega_2}\tilde\phi_k\geq \hat\phi_k\geq\rho'-k$$ and therefore the function $v$ given by $$v(z)=\left\{
\begin{array}{ll}
    P_{\Omega_2}\tilde\phi_k(z)    & \hbox{ for  } z\in\Omega_2\\
    \phi_k(z)=\rho'(z)-k                & \hbox{ on }  \Omega\setminus\Omega_2,
\end{array}\right. $$ is plurisubharmonic on $\Omega$ and smaller than $\tilde{\phi_k}$. Thus $ P_{\Omega_2}\tilde\phi_k=\hat\phi_k|_{\Omega_2}$. Observe that for $z\in\Omega_2$ and $w\in K_\varepsilon$ we have $\tilde{\phi}_k(z)\geq\tilde{\phi}_k(z+w)-\omega(|w|)$, where $\omega$ is the modulus of continuity of a function $\tilde{\phi}_k|_L$ and $L=L_\varepsilon+K_\varepsilon=\{z+w:z\in L_\varepsilon \hbox{ and } w\in K_\varepsilon\}$. Therefore, $\hat\phi(z)\geq\hat\phi(z+w)-\omega(|w|)$. By Lemma \ref{ciaglosc} the function $\hat\phi_k|_{\Omega_2}$ is continuous.
 
Because any $z\in\Omega$ is in $\Omega_2(\varepsilon)$ for some $\varepsilon$ as above, we obtain that the function $\hat\phi_k$ is continuous on $\Omega$. Using the Richberg theorem we can modify the sequence $\hat\phi_k$ to a sequence $u_k$ of smooth plurisubharmonic functions which decreases to $u$. $\Box$

The approximation by continuous functions, can be proved in the same way in a much more general situation.

 \begin{thr}\label{apFsub} Let $\textbf{F}$ be a constant coefficient subequetion such that all $\textbf{F}$-subharmonic functions are subharmonic and all convex functions are $\textbf{F}$-subharmonic. Let $D\subset\mathbb{R}^n$ be a domain with  $Lipschitz$ boundary. Then for any  point $P\in \bar D$ there is a neighbourhood $U$ of $P$ such that for any function $u\in \textbf{F}(D\cap U)$ there is a sequence $(u_k)\subset\textbf{F}(D\cap U)$ of continuous functions decreasing to $u$.
\end{thr}

Here we use terminology from \cite{h-l}\footnote{See also theorem 2.6 there for elementary properties of $\textbf{F}$-subharmonic functions needed in the proof. For the result of local continuous approximation of $F$-subharmonic functions see \cite{h-l-p}}.

\section{Example}\label{Example}

Similarly as in Lecture 14 in \cite{f-s} we construct a  domain $\Omega\subset\cn$ and a plurisubharmonic function $u$ on $\Omega$ which can not be regularise. Let
$A=\{\frac{1}{k}:k\in\mathbb{N}\}$ and a sequence $(x_k)\subset(0,1)\setminus A$ is such that its limit set is equals $\bar A$. Put
$$\lambda(z)=\sum_{k=1}^\infty c_k\log|z-x_k| \hbox{ for } z\in\mathbb{C},$$
where $c_k$ is a such sequence of numbers rapidly decreasing to 0, such that\\
i) $\lambda$ is a subharmonic function on $\mathbb{C}$ and\\
ii)$\lambda|_A\geq-\frac{1}{2}$.\\
For $k\in\mathbb{N}$ let $D_k$ be a disc with  center  $x_k$ such that $\lambda|_{D_k}<-1$. Now, we can define
$$\Omega=\{(z',z_n)\in\mathbb{C}^{n-1}\times\mathbb{C}:|z'|^2+|z_n-1|^2<1\}\setminus K,$$ 
where
$$K=\{(z',z_n)\in\mathbb{C}^{n-1}\times\mathbb{C}:|z'|=|z_n| \hbox{ and } z_n\notin\cup_{k=1}^\infty D_k\},$$
and
$$u(z',z_n)=\left\{
\begin{array}{ll}
    -1    & \hbox{ for  } z\in\Omega\cap D\\
    \max\{\lambda{(z_n)},-1\}                & \hbox{ on }  \Omega\setminus D,
\end{array}\right.$$
where $D=\{(z',z_n)\in\mathbb{C}^{n-1}\times\mathbb{C}:|z'|<|z_n| \}$.

Let $U$ be any neighbourhood of $0$. We can choose numbers $k\in\mathbb{N}$, $0<r<\frac{1}{k}$ such that $$\Omega_U:=\{(z',z_n)\in\mathbb{C}^{n-1}\times\mathbb{C}:|z'|\leq\frac{2}{k},|z_n-\frac{1}{k}|\leq r\}\setminus K\subset \Omega\cap U.$$
 Let $y_p$ be subsequence of $x_p$ such that $\lim_{p\rightarrow\infty}y_p=\frac{1}{k}$ and for all $p$ we have $|y_p-\frac{1}{k}|<r$. Now, we can repeat  the argument from  \cite{f-s}. If $u_q$ is a sequence of smooth plurisubharmonic functions decreasing to $u$ on $\Omega_U$, then for $q$ sufficiently large $u_q\leq-\frac{3}{4}$ on the set
$$\{(z',z_n)\in\mathbb{C}^{n-1}\times\mathbb{C}:|z'|=\frac{2}{k},|z_n-\frac{1}{k}|\leq r\}\subset\partial\Omega_U.$$
By the maximum principle (on  sets $\{(z',z_n)\in\mathbb{C}^{n-1}\times\mathbb{C}:|z'|\leq\frac{2}{k},z_n=y_p\}\subset\mathbb{C}^{n-1\times}\{y_p\}$) we also have $u_q(0,y_p)\leq-\frac{3}{4}$ and by continuity of $u_q$ we get $u_q(0,\frac{1}{k})\leq-\frac{3}{4}<u(0,\frac{1}{k})$. This is a contradiction.

Note that we have not only just proved that for the above $\Omega$ Theorem \ref{aproksymacjanabrzegu} does not hold but we also have the following stronger result:
\begin{prp}
 Let $\Omega$ and $u$ be as above. Then  the function $u$  can not be smoothed on $U\cap\Omega$ for any neighbourhood $U$ of $0\in\partial\Omega$.
\end{prp}

\section{Questions}
In this last section we state some open questions related to the content of the note.

{\bf  1.} Is it possible to characterize to characterize $\mathcal{S}$-domains? In view of Corollary \ref{lok}, even in the class of smooth domains, it seems to be a challenging problem.  

{\bf 2.} Let $D$ be  an $\mathcal{S}$-domain and let $f$ be a continuous function which is bounded from below. Is the plurisubharmonic envelope of $f$ continuous? Assume in addition that $f$ is smooth.   What is the optimal regularity of $P_Df$? The author does not know answers even in the case of the ball in $\cn$. Note  that if  $D$ is not an $\mathcal{S}$-domain, then there is a smooth function $f$ bounded from below such that  $P_Df$ is discontinuous.

{\bf 3.} What is the optimal assumption about the regularity of the boundary of $D$ in the Theorem \ref{aproksymacjanabrzegu}? Is it enough to assume that $D=int\bar D$?

{\bf 4.} Let $M$ be a real (smooth or Lipschitz) hypersurface in $\cn$.    Is it true that for any $P\in M$ there exists a smooth pseudoconvex neighbourhood $U\subset B$ such that $M$ divides $U$ into two $\mathcal{S}$-domains?

\end{document}